\documentclass[10pt,twoside]{article}
\usepackage{amssymb,amscd}
\usepackage{graphicx}
\usepackage{amsmath}
\usepackage{Latex-document}

\def\volumeno{I}

\title{\bf Hidden Markov and State Space Models \vskip -2mm Asymptotic Analysis of
Exact and \vskip -2mm Approximate Methods for Prediction, \vskip -2mm
Filtering, Smoothing and Statistical \vskip -2mm Inference\vskip 6mm}

\author{{\bf P. Bickel}\thanks{University of California, Berkeley, USA. E-mail: bickel@stat.berkeley.edu}
\quad Y. Ritov\thanks{Hebrew University, Israel} \quad T.
Ryden\thanks{University of Lund, Sweden}\vspace*{-0.5cm}}
\date{\vspace{-8mm}}

\begin{document}

\maketitle

\thispagestyle{first} \setcounter{page}{555}

\begin{abstract}

\vskip 3mm

State space and hidden Markov models can both be subsumed under the same
mathematical structure.  On a suitable probability space
$(\Omega,{\mathcal A},P)$ are defined
$(X_1,Y_1,X_2,Y_2,\dots,X_n,Y_n,\dots)$ a sequence of random
``variables'' taking values in a product space $\prod_{j=1}^{\infty}
({\mathcal X}_j \times {\mathcal Y}_j)$ with an appropriate sigma field.
The joint behavior under $P$ is that the $X_j$ are stationary Markovian
and that given $(X_1,X_2,\dots)$ the $Y_j$ are independent and further
that $Y_j$ is independent of all $X_i: i \ne j$ given $X_j$. If
${\mathcal H}$ is finite these are referred to as Hidden Markov models.
The general case though focussing on ${\mathcal X}$ Euclidean is
referred to as state space models.  Essentially we observe only the
$Y$'s and want to infer statistical properties of the $X$'s given the
$Y$'s.  The fundamental problems of filtering, smoothing prediction are
to give algorithms for computing exactly or approximately the
conditional distribution of $X_t$ given $(Y_1,\dots,Y_t)$ (Filtering),
the conditional distribution of $X_t$ given $Y_1,\dots,Y_T$, $T > t$
(Smoothing) and the conditional distribution of $X_{t+1},\dots,X_T$
given $Y_1,\dots,Y_t$ (Prediction).  If as is usually the case $P$ is
unknown and is assumed to belong to a smooth parametric family of
probabilities $\{P_{\theta}: \theta \in R^d\}$ ,we face the further
problem of efficiently estimating $\theta$ using $Y_1,\dots,Y_T$
(computation of the likelihood, and maximum likelihood estimation,
etc.).

State space models have long played an important role in signal
processing.  The Gaussian case can be treated algorithmically using the
famous Kalman filter [6]. Similarly since the 1970s there has been
extensive application of Hidden Markov models in speech recognition with
prediction being the most important goal.  The basic theoretical work
here, in the case ${\mathcal X}$ and ${\mathcal Y}$ finite (small)
providing both algorithms and asymptotic analysis for inference is that
of Baum and colleagues [1].  During the last 30-40 years these general
models have proved of great value in applications ranging from genomics
to finance---see for example [7].

Unless the $X,Y$ are jointly Gaussian or ${\mathcal X}$ is finite and
small the problem of calculating the distributions discussed and the
likelihood exactly are numerically intractable and if ${\mathcal Y}$ is
not finite asymptotic analysis becomes much more difficult.  Some new
developments have been the construction of so-called ``particle
filters'' (Monte Carlo type) methods for approximate calculation of
these distributions (see Doucet et  al. [4]) for instance and general
asymptotic methods for analysis of statistical methods in HMM [2] and
other authors.

We will discuss these methods and results in the light of exponential
mixing properties of the conditional (posterior) distribution of
$(X_1,X_2,\dots)$ given $(Y_1,Y_2,\dots)$ already noted by Baum and
Petrie [1] and recent work of the authors Bickel, Ritov and Ryden [3],
Del Moral and Jacod in [4], Douc and Matias [5].

\vskip 4.5mm


%\noindent {\bf Keywords and Phrases:} Mirror symmetry, D-branes, Floer
%homology.
\end{abstract}

\label{lastpage}

\end{document}